\documentclass[11pt]{article}
\usepackage[english]{babel}
\usepackage[T1]{fontenc}
\usepackage{times}
\usepackage{graphicx}
\usepackage{epsfig}
\usepackage{amsfonts,amsmath}
\usepackage{latexsym}
\usepackage{amssymb}
\usepackage{theorem}
\usepackage{color}
\usepackage{hyperref}
\def \qed{\hspace*{\fill} $\Box$\par\medskip}

\def \xtx {X^{t,x}}

\def \vy {\vec y}

\def \mi {\mathcal{I}}

\def \txk {(t,x) \in [0,T]\times \R^{k}}
\def \spa { [0,T]\times \R^{k}}
\def \spat { [t,T]\times \R^{k}}

\def \hdd{{\cal H}^{2,d}}
\def \d {\delta}
\def \wt {\mbox{w.r.t}}

\def \ss {{\cal S}^2}
\def \ii {i\in \mi}
\def \lb {\label}
\def \tx {(t,x)\in [0,T]\times \R^k}
\def \sol{\underline}
\def \nd {\noindent}
\def \tbf{\textbf}

\newlength{\oldparindent}
\makeatletter
\newcommand*{\rom}[2]{\expandafter\@slowromancap\romannumeral #1@}
\newcommand{\E}{\mathbb{E}}

\newcommand{\R}{\mathbb{R}}
\newcommand{\N}{\mathbb{N}}

\def \ssp {\mbox{esssup}}
\numberwithin{equation}{section}
\newtheorem{definition}{Definition }[section]
\newtheorem{proposition}[definition]{Proposition }
{Lemma }
\newtheorem{theoreme}[definition]%
{Theorem }
{Corollary }
\newtheorem{remarque}[definition]%
{Remark }
{Assumption}
\newcommand{\MBFigure}[6]{
$\left. \right.$ \\
\refstepcounter{figure}
\addcontentsline{lof}{figure}{\numberline{\thefigure}{\ignorespaces #5}}
\begin{center}
\begin{minipage}{#1cm}
\centerline{\includegraphics[width=#2cm,angle=#3]{#4}}
\begin{center}
\upshape{F\textsc{ig} \normal
\end{center}
size{\thefigure}. $-$} #5
\end{center}
\label{#6}
\end{minipage}
\end{center}
$\left. \right.$ \\}

\title{\Large \bf Viscosity Solutions of Systems of PDEs with Interconnected
Obstacles and  Switching Problem without Monotonicity Condition}
\author{{\Large Said Hamad\`ene}\thanks{LMM, Le Mans Universit\'e, Avenue Olivier Messiaen, 72085 Le Mans, Cedex 9, France. \texttt{e-mail: hamadene@univ-lemans.fr};} {,\Large }
 {\,\,\Large Mohamed Mnif}\thanks{ University of Tunis El Manar, Laboratoire de Mod\'elisation Math\'ematique et Num\'erique  e-mail: \texttt{mohamed.mnif@enit.rnu.tn};}
 {\,\,\Large and }
 {\Large Sarra Neffati.}\thanks{University of Tunis El Manar,
 Laboratoire de Mod\'elisation Math\'ematique et Num\'erique e-mail: \texttt{sarra.neffati@enit.utm.tn}.}}
\begin{document}
\maketitle
\begin{abstract}
We show the existence and uniqueness of a continuous  viscosity solution of a system of
partial differential equations (PDEs for short) without assuming the usual monotonicity
conditions on the driver function as in Hamad\`ene and Morlais's article \cite{hamadene2013viscosity}. Our method strongly relies on the link between PDEs and  reflected backward stochastic differential equations with interconnected obstacles for which we already know that the solution exists and is unique for general drivers.
\end{abstract}

\textbf{Keywords}: Partial differential equations ; Interconnected obstacles ; Viscosity solution ; Multi-modes switching ; HJB system ; Reflected Backward stochastic differential equations.

\section{Introduction}

 \makebox[0.5cm][r] The main objective of this paper is to study the problem of existence and uniqueness of a solution in viscosity sense $(u^i)_{i=1,m}$ of the following system of partial differential equations with obstacles which depend on the solution:  $\forall i \in \mathcal{I}:= \lbrace 1,...,m\rbrace,$
\begin{equation}\label{eqq1}
\begin{cases}
\min \lbrace u^i(t,x) - \displaystyle \max_{j \in {\mathcal{I}^{-i}}}\{u^j(t,x)-g_{ij}(t,x)\} ;\\
 \quad \quad -\partial_tu^i(t,x) - \mathcal{L}u^i(t,x) - f_i(t,x,(u^k(t,x))_{k=1,...,m},(\sigma^\top D_xu^i)(t,x))\rbrace = 0 \,\, ;\\
u^i(T,x) = h_i(x)
\end{cases}
\end{equation}
where $\mathcal{I}^{-i}:= \mathcal{I}-\lbrace i \rbrace$ and $\mathcal{L}$ is an infinitesimal generator which has the following form
\begin{equation}\label{generinfintro}
\mathcal{L}\varphi(t,x):= b(t,x)^\top.D_x\varphi(t,x) + \frac{1}{2}\mbox{Tr}[(\sigma\sigma^\top)(t,x)D_{xx}^2\varphi(t,x)]
\end{equation}and which is associated with the stochastic process $X^{t,x}$ solution of the SDE \eqref{xtx}. 

As pointed out previously, in \eqref{eqq1}, the obstable of $u^i$ is the function 
$\max_{j \in {\mathcal{I}^{-i}}}\{u^j(t,x)-g_{ij}(t,x)\}$ which actually depends on 
the solution $(u^i)_{i=1,m}$, which means that the obstacles are interconnected. 

This problem is related to the optimal stochastic switching control problem which can be described, through an example, as follows: Assume that we have a power plant which has several modes of production and which the manager puts in a specific mode according to its profitability which depends on the electricity price in the energy market evolving according to the following stochastic differential equation
\begin{equation}\label{xtx}
dX_s^{t,x} = b(s, X_s^{t,x}) ds + \sigma(s, X_s^{t,x}) dB_s, \hspace{0.2cm}  s\geq t \mbox{ and }X^{t,x}_t=x .
\end{equation}
The aim of the manager is to maximize her global profit over an horizon $[0,T]$ by optimally choosing controls of the form $\delta := (\theta_k, \alpha_k)_{k \geqslant 0}$ where $(\theta_k)_{k \geqslant 0}$ is an increasing sequence of stopping times  at which the manager  switches the system  across the different operating modes and $(\alpha_k)_{k \geqslant 0}$ is a sequence of random variables with values in $\lbrace 1,...,m\rbrace$ which stand for the modes to which the production is switched. Namely for any $k\geq 1$, at $\theta_k$, the manager switches the production from $\theta_{k-1}$ to $\theta_{k}$ ($\theta_0$ and $\alpha_0$ are the starting time and mode respectively). However, switching the plant from the mode $i$ to the mode $j$ is not free generates expenditures which amount to $g_{ij}(s,X_s^{t,x})$ at time $s$. When the plant is run under a strategy $\delta$, its yield is given by
\begin{equation*}\begin{array}{c}
J({\delta};t,x) := \mathbb{E}\Big[ \int_t^T
f^{\delta}(s,X_s^{t,x})ds   - A_T^{\delta} +
h^{\delta}(X_T^{t,x})\Big]\end{array}
\end{equation*}
where:

\nd a) $f^{\delta}(s,X_s^{t,x})$ is the instantaneous payoff of the station when run under $\d$ and $h^{\delta}(X_T^{t,x})$ is the terminal payoff ;

\nd b) the quantity $A_T^{\delta}$ stands for the total switching cost when the strategy $\delta$ is implemented (see \eqref{defcoutswitch} for its definition).

\noindent The problem is to find an optimal management strategy $\delta^*$, i.e.,
which satisfies \\ $J({\delta}^*;t,x) = \sup \{J({\delta};t,x), \delta\in {\cal A}_{ad}\} $. This latter quantity is nothing but the fair price of the power station in the energy market.

In \eqref{eqq1}, if for any $\ii$, $f_i$ does not depend on
$(u^k)_{k=1,m}$ and $D_xu^i$, the system reduces to the
Hamilton-Jacobi-Bellman one associated with the switching problem
and it is shown in \cite{hamadene2013viscosity} that it has a
unique solution $(u^i)_{i=1,m}$ which satisfies 
$$
u^i(t,x)=\sup \{J({\delta};t,x), \delta\in {\cal A}_t^i\}$$ where
${\cal A}_t^i$ is the set of admissible strategies which start at
time $t$ from mode $i$. 

In a so general form, system \eqref{eqq1} can be related 
to switching problems with utility functions \cite{hamadene2009optimal}, recursive utilities, knightian uncertainty \cite{hamadene2010switching}, etc.   

The main tool to tackle system \eqref{eqq1} is to deal with the
following system of reflected backward stochastic differential
equations (RBSDEs for short) with interconnected obstacles:
$\forall i \in \lbrace 1,...,m\rbrace$ and $s \in [t,T],$
\begin{equation}\label{eq17intro}
\begin{cases}
Y_s^{i,t,x} = h_i(X_T^{t,x}) + \int_s^T f_i(r,X_r^{t,x},(Y_r^{k,t,x})_{k=1,,...,m},Z_r^{i,t,x})dr + K_T^{i,t,x} - K_s^{i,t,x}- \int_s^T Z_r^{i,t,x}dB_r, \\
  Y_s^{i,t,x} \geqslant  \displaystyle \max_{j\in \mathcal{I}^{-i}}(Y_s^{j,t,x} -g_{ij}(s,X_s^{t,x})), \\
 \int_t^T [Y_s^{i,t,x} -\max_{j\in \mathcal{I}^{-i}}( Y_s^{j,t,x} -g_{ij}(s,X_s^{t,x}))] dK_s^{i,t,x} = 0.
\end{cases}
\end{equation}
 This system of RBSDEs has been investigated in several papers including
 (\cite{chassagneux2011note,hamadene2010switching,hamadene2013viscosity, hu2010multi}, etc.). In \cite{chassagneux2011note}, the authors proved
 that it has a unique solution
 $(Y^{i,t,x},Z^{i,t,x},K^{i,t,x})_{\ii}$ if the functions $(f_i)_{\ii}$ are
 merely Lipschitz w.r.t $((y_l)_{l=1,m},z)$.

 Concerning now the system of PDEs \eqref{eqq1}, Hamadene et al. proved in
 \cite{hamadene2013viscosity}, that
 if for any $\ii$ and $k\in {\cal I}^{-i}$, $f_i(t,x,(y_l)_{l=1,m},z)$ is increasing w.r.t
 $y_k$ (see (H4)-(i)), then system \eqref{eq17intro} has a unique solution $(u^i)_{i=1,m}$ in the class
 of continuous functions with polynomial growth and which is given by:
 \begin{equation}\label{fkc}\forall \ii, \,\, u^{i}(t,x)=
 Y^{i;t,x}_t,\,\,\tx,\end{equation} where $(Y^{i;t,x})_{\ii}$ is the first
 component of the solution of the system of reflected BSDEs
 \eqref{eq17intro}. The same result is obtained if, instead of $(f_i)_{\ii}$, their opposites
 $(-f_i)_{\ii}$ verify the previous monotonicity property (see (H4)-(ii)). However
 without assuming those monotonicity conditions on the drivers $(f_i)_{\ii}$ the
 problem of existence and uniqueness of the solution in viscosity
 sense of system \eqref{eq17intro} remains open. In this paper we show that system \eqref{eq17intro}
 has a unique solution
 without assuming the previous monotonicity properties on the drivers $(f_i)_{\ii}$. This is the main novelty of this work. As a consequence,
 we make matching the probabilistic and PDEs frameworks. Once more
 our method relies on the link between reflected BSDEs and PDEs with
 obstacles in the Markovian framework of randomness.

The paper is organized as follows. In Section 2, we formulate
accurately the problem. In section 3,  we show that Feynman-Kac
formula holds for the components $(Y^{i;t,x})_{\ii}$ of the solution
of \eqref{eq17intro}, i.e., the representation \eqref{fkc} holds true.
In Section 4, we show that the functions $(u^i)_{\ii}$ are
continuous and are the unique viscosity solution of \eqref{eqq1} in
the class of functions with polynomial growth. The proof is deeply
related to the fact that system \eqref{eq17intro} of RBSDEs has a
unique solution. \qed
\section{Preliminaries and notations}
Let $T > 0$ be a given time horizon and $(\Omega, \mathcal{F}, \mathbb{P})$  be a probability space on which is defined a standard $d$-dimensional Brownian motion $B = (B_t)_{t\leq T}$ whose natural filtration is $(\mathcal{F}^0_t:= \sigma(B_s, s\leqslant t))_{t\le T}$ and $\mathbb{F}= (\mathcal{F}_t)_{0\leq t \leq T}$ is its augmentation with the $\mathbb{P}$-null sets of $\mathcal{F}$. Hence $(\mathcal{F}_t)_{0\leq t \leq T}$ is right continuous and complete.

 We now introduce the following spaces :
 \begin{itemize}
 \item[a)] $\mathcal{P}$  is the $\sigma$-algebra of $\mathbb{F}$-progressively measurable sets on
 $ [0,T]\times \Omega ;$
\item[b)] $\mathcal{S}^2$ is the set of $\mathcal{P}$-measurable, continuous, $\R$-valued  processes $Y=(Y_s)_{s\leq T}$ such that $\mathbb{E}[\displaystyle \sup_{s\leq T} |Y_s|^2] < \infty ;$
  \item[c)] $\mathcal{A}^2 $ is the subset of $\mathcal{S}^2$ of non decreasing processes $K=(K_t)_{t\leq T}$ such that $K_0 = 0 $ ;
\item[d)] $\mathcal{H}^{2,l}$ $(l\geqslant 1)$ is the set of $\mathcal{P}$-measurable and $\R^l$-valued processes $Z:=(Z_s)_{s\leq T}$ such that $\mathbb{E}[\int_0^T |Z_s|^2 ds] < \infty$.
 \end{itemize}

  Next, for any given $(t,x) \in [0,T]\times \R^k$ ($k$ is a positive integer), we consider the following standard stochastic differential equation (SDE) :
\begin{equation}\label{eq1}
 \begin{cases}
dX_s^{t,x} = b(s, X_s^{t,x}) ds + \sigma(s, X_s^{t,x}) dB_s, \hspace{0.2cm}  s\in [t, T] \\
 X_s^{t,x} = x, \hspace{0.2cm} 0\leq s\leq t
 \end{cases}
 \end{equation}
 where $b : [0,T]\times \R^k \rightarrow \R^k$ and $\sigma : [0,T]\times \R^k \rightarrow \R^{k\times d}$ are two continuous functions and Lipschitz  $\wt$ $x$, i.e., there exists a positive constant $C$ such that
 \begin{equation}\label{eq2}
 |b(t,x) - b(t,x^{\prime})| + |\sigma(t,x) - \sigma(t,x^{\prime})|\leq C|x-x^{\prime}|,\hspace{0.2cm} \forall (t,x,x^{\prime}) \in [0,T]\times \R^{k+k}.
 \end{equation}
 Note that the continuity of $b$, $\sigma$ and \eqref{eq2} imply the existence of a constant $C$ such that  \begin{equation}\label{eq3}
 |b(t,x)| + |\sigma(t,x)| \leq C(1+|x|), \hspace{0.2cm} \forall (t,x) \in [0,T]\times \R^{k}.
\end{equation}
Conditions \eqref{eq2} and \eqref{eq3} ensure, for any $\txk$, the existence and uniqueness of a solution $\lbrace X_s^{t,x}, t\leq s\leq T\rbrace$ to the SDE \eqref{eq1} (see \cite{revuz1999continuous} for more details). Moreover, it satisfies the following estimate: $\forall p\geqslant 1,$\\
\begin{equation}
\mathbb{E}[\sup_{s\leq T} |X_s^{t,x}|^p] \leqslant C(1+|x|^p).
\end{equation}

Next let us introduce the following deterministic functions $(f_i)_{i=1,...,m}$, $(h_i)_{i=1,...,m}$ and $(g_{ij})_{i,j=1,...,m}$ defined as follows : for any $i,j \in \lbrace 1,...,m \rbrace$,
\begin{align*}
&a)\,\,\,f_i : (t,x,\vec{y},z) \in [0,T] \times \R^{k+m+d}  \longmapsto f_i(t,x,\vec y,z)\in \R \,\, (\vec y:=(y^1,...,y^m));\\\\
&b)\,\, g_{ij}: \txk \longmapsto g_{ij}(t,x)\in \R \quad ; \\\\&c)\,\,h_i :x\in \R^{k}  \longmapsto h_i(x)\in \R.
\end{align*}
Additionally we assume that they satisfy:
\begin{itemize}
\item[\textbf{(H1)}] For any $i \in \lbrace 1,...,m\rbrace$,
\begin{itemize}
 \item[(i)]  The function $(t,x) \mapsto f_i(t,x,\vy ,z)$  is continuous, uniformly w.r.t. the variables $(\vy, z)$,
 \item[(ii)] The function $f_i$ is Lipschitz continuous with respect to the variables $(\vy,z)$ uniformly in $(t,x)$, i.e., there exists a positive constant $C_i$ such that for any $(t,x) \in [0,T]\times \R^k,$
 $ (\vy, z)$ and $(\vy_1,z_1)$ elements of  $\R^{m+d}$:
\begin{equation}
|f_i(t,x,\vy,z)- f_i(t,x,\vy_1,z_1)| \leq C_i( |\vy- \vy_1| + |z - z_1|).
\end{equation}
\item[(iii)] The mapping $(t,x) \mapsto f_i(t,x,0,...,0)$ has polynomial growth in $x$, i.e., there exist two constants $C > 0$ and $p \geqslant 1$ such that for any $(t,x) \in [0,T]\times \R^k$,
\begin{equation}
|f_i(t,x,0,...,0)| \leq C(1+|x|^p).
\end{equation}
\end{itemize}
\item[\textbf{(H2)}] $\forall i,j \in \lbrace 1,...,m \rbrace,$ $g_{ii} = 0$ and for $i\neq j,$ $g_{ij}(t,x)$ is non-negative, jointly continuous in $(t,x)$ with polynomial growth and satisfy the following non free loop property :

For any $(t,x) \in [0,T] \times \R^k$, for any sequence of indices $i_1,...,i_k$ such that $i_1 = i_k$ and $card\lbrace i_1,...,i_k \rbrace = k-1$ ($k\ge 3$) we have
\begin{equation}
g_{i_{1}i_{2}}(t,x) + g_{i_{2}i_{3}}(t,x) + ... + g_{i_{k}i_{1}}(t,x)> 0.
\end{equation}
 \item[\textbf{(H3)}] For $i \in \lbrace 1,...,m\rbrace$, the function $h_i$, which stands for the terminal condition, is continuous with polynomial growth and satisfies the following consistency condition:
 \begin{equation}
\forall x \in \R^k,\,\, h_i(x) \geqslant \max_{j\in \mathcal{I}^{-i}}\{h_j(x)-g_{ij}(T,x)\}.
 \end{equation}

\item[\textbf{(H4)-(i)}] $\forall i \in \mathcal{I}\mbox{ and } j\in \mathcal{I}^{-i}$, the mapping $
w \in \R \longmapsto f_i(t,x,y^1,...,y^{j-1},w,y^{j+1},...,y^m,z)$ is non-decreasing whenever the other components   $(t,x,y^1,...,y^{j-1},y^{j+1},...,y^m,z)$ are fixed.
\item[\textbf{(H4)-(ii)}] the functions $(-f_i)_{\ii}$ verify (H4)-(i). \qed
\end{itemize}

The main objective of this paper is to study the following system of PDEs with interconnected obstacles: For any $i \in \mathcal{I}:=\lbrace 1,...,m\rbrace$,
\begin{equation}\label{eqq12}
\begin{cases}
\min \lbrace u^i(t,x) - \displaystyle \max_{j \in {\mathcal{I}^{-i}}}(u^j(t,x)-g_{ij}(t,x)) ;\\
 \quad \quad -\partial_tu^i(t,x) - \mathcal{L}u^i(t,x) - f_i(t,x,(u^k(t,x))_{k=1,...,m},(\sigma^\top D_xu^i)(t,x))\rbrace = 0 \,\, ;\\
u^i(T,x) = h_i(x)
\end{cases}
\end{equation}where the operator $\mathcal{L}$ is the infinitesimal generator associated with $\xtx$, i.e., \begin{equation}\label{generinf}
\mathcal{L}\varphi(t,x):= b(t,x)^\top.D_x\varphi(t,x) + \frac{1}{2}\mbox{Tr}[(\sigma\sigma^\top)(t,x)D_{xx}^2\varphi(t,x)]
\end{equation}for any $\R$-valued function $\varphi(t,x)$ such that $D_x\varphi(t,x)$ and $D^2_{xx}\varphi(t,x)$ are defined.

A solution $(u^i)_{i\in \mi}$ of system \eqref{eqq12} is to be understood in viscosity sense whose definition is the following:
\begin{definition}
Let $\vec u:=(u^i)_{i\in \mi}$ be a function of $C([0,T] \times \R^k ; \R^m)$. We say that $\vec u$ is a viscosity supersolution (resp. subsolution) of \eqref{eqq1} if:
$\forall i \in \lbrace 1,...,m\rbrace$,
\begin{align*}
& \mbox{a) } \,u^i(T,x) \geq (\mbox{resp.} \leq )\,\, h_i(x) , \,\,\forall x \in \R^k\,\,; \\
& \mbox{b)  if}\,\phi \in {\cal C}^{1,2}([0,T] \times \R^k)\mbox{  is
such that $(t,x) \in [0,T) \times \R^k$ is a local minimum (resp. maximum) point of $u^i - \phi$ then }
\end{align*}
\begin{align*}\label{systemdef}
\min \Big\lbrace &u^i(t,x) - \displaystyle \max_{j \in {\mathcal{I}^{-i}}}(u^j(t,x)-g_{ij}(t,x)) ;\\
 &-\partial_t\phi(t,x) - \mathcal{L}\phi(t,x) - f_i(t,x,(u^k(t,x))_{k=1,...,m},(\sigma^\top D_x\phi)(t,x))\Big\rbrace \geq \,\,(resp. \le )\,\,0.
\end{align*}
(ii) We say that $\vec u:=(u^i)_{i\in \mi}$ is a viscosity solution of \eqref{eqq1} if it is both  a supersolution and subsolution of \eqref{eqq1}.
\end{definition}
 \section{Connection with Systems of Reflected BSDEs with Oblique Reflection}
 The viscosity solution of system \eqref{eqq12} is deeply connected (one can see \cite{hamadene2013viscosity} for more details) with the following system of reflected BSDEs with interconnected obstacles (or oblique reflection) associated with \\$((f_i)_{i\in \mathcal{I}}, (g_{ij})_{i,j \in \mathcal{I}},(h_i)_{i \in \mathcal{I}})$ : $\forall i = 1,...,m$ and $s \in [t,T]$,
\begin{equation}\label{eqq2}
\begin{cases}
\textstyle{
Y^{i,t,x} \in \mathcal{S}^2,\hspace{0.1cm} Z^{i,t,x} \in \mathcal{H}^{2,d}} \hspace{0.1cm} $and$ \hspace{0.1cm}  \textstyle {K^{i,t,x} \in \mathcal{A}^2;} \\\\
\textstyle {Y_s^{i,t,x} = h_i(X_T^{t,x}) + \int_s^T f_i(r,X_r^{t,x},(Y_r^{k,t,x})_{k=1,...m},Z_r^{i,t,x})dr + K_T^{i,t,x} - K_s^{i,t,x}- \int_s^T Z_r^{i,t,x}dB_r,}\\\\
 \textstyle {Y_s^{i,t,x} \geqslant  \displaystyle \max_{j \in \mathcal{I}^{-i}}(Y_s^{j,t,x} -g_{ij}(s,X_s^{t,x}))}\mbox{ and }\textstyle {\int_t^T (Y_s^{i,t,x} -\displaystyle \max_{j\in \mathcal{I}^{-i}}(Y_s^{j,t,x} -g_{ij}(s,X_s^{t,x}))) dK_s^{i,t,x} = 0.}
\end{cases}
\end{equation}This system \eqref{eqq2} of reflected BSDEs is considered in several works (see e.g. \cite{hamadene2010switching,hamadene2013viscosity,hu2010multi,chassagneux2011note}, etc.). Under (H1)-(H3) and (H4)-(i) as well, this system has been considered first in \cite{hamadene2010switching} where issues of existence and uniqueness of the solution, and comparison of the solutions, are considered. Actually it is shown:
\begin{theoreme} \label{existmono}(see \cite{hamadene2010switching}).\\
i) Assume that the deterministic functions $(f_i)_{i\in \mathcal{I}}, (g_{ij})_{i,j \in \mathcal{I}}$ and $(h_i)_{i \in \mathcal{I}}$ verify Assumptions (H1)-(H3) and (H4)-(i). Then system \eqref{eq12} has a unique solution
$(Y^i,Z^i,K^i)_{i\in \mi}$.\\
ii) If
$(\bar f_i)_{i\in \mathcal{I}}, (\bar g_{ij})_{i,j \in \mathcal{I}}$ and $(\bar h_i)_{i \in \mathcal{I}})$ are other functions satisfying (H1)-(H3) and (H4)-(i) and, moreover, for any $i,j\in \mi$,
$$
f_i\leq \bar f_i, \,\,h_i\leq \bar h_i \mbox{ and }g_{ij}\ge \bar g_{ij}. $$
Then for any $i\in \mi$, $Y^i\leq \bar Y^i$ where $(\bar Y^i,\bar Z^i,\bar K^i)_{i\in \mi}$ is the solution of the system associated with $(\bar f_i)_{i\in \mathcal{I}}, (\bar g_{ij})_{i,j \in \mathcal{I}}$ and $(\bar h_i)_{i \in \mathcal{I}}$. \qed
\end{theoreme}
In \cite{chassagneux2011note}, Chassagneux et al. have also considered system \eqref{eqq2} without assuming Assumption (H4)-(i). They stated the following result:
\begin{theoreme}\label{thmexistence} (see \cite{chassagneux2011note}) Assume that the deterministic functions $(f_i)_{i\in \mathcal{I}}, (g_{ij})_{i,j \in \mathcal{I}}$ and $(h_i)_{i \in \mathcal{I}}$ verify Assumptions (H1)-(H3). Then system \eqref{eq12} has a unique solution
$(Y^i,Z^i,K^i)_{i\in \mi}$.
\end{theoreme}
{\bf \sol{Proof}}: We give the main steps of the proof of this result since it plays an important role in the proof of our main result. This proof is mainly based on the interpretation of the solutions of \eqref{eqq2},
when $(f_i)_{i=1,m}$ do not depend on $\vec y$, as the value function of an optimal switching problem. Indeed let $\vec{\Gamma} := (\Gamma^i)_{i=1,...,m} \in \mathcal{H}^{2,m}$ and let us introduce the following mapping:
\begin{align}\label{defteta}
\Theta : \hspace{0.1cm} & \mathcal{H}^{2,m} \rightarrow \mathcal{H}^{2,m} \cr
&\vec{\Gamma} \quad \mapsto \Theta(\vec{\Gamma}):=(Y^{\Gamma,i})_{i=1,...,m}
\end{align}
  where $(Y^{\Gamma,i},Z^{\Gamma,i},K^{\Gamma,i})_{i \in \mathcal{I}} \in (\mathcal{S}^2 \times \mathcal{H}^{2,d} \times \mathcal{A}^2)^m$ (we omit the dependence on $t,x$ of
$Y^{\Gamma,i},Z^{\Gamma,i},K^{\Gamma,i}$ as no confusion is possible) is the unique solution of the following system of reflected BSDEs with interconnected obstacles (or oblique reflection): $\forall  i \in \mathcal{I}$,
 \begin{equation} \label{eq12}
\begin{cases}
\textstyle {Y_s^{\Gamma,i} = h_i(X_T^{t,x}) + \int_s^T f_i(r,X_r^{t,x},\vec{\Gamma}_r,Z_r^{\Gamma,i})dr + K_T^{\Gamma,i} - K_s^{\Gamma,i}- \int_s^T Z_r^{\Gamma,i}dB_r, \quad\forall s \leq T;}\\\\
 \textstyle {Y_s^{\Gamma,i} \geqslant  \displaystyle \max_{j \in \mathcal{I}^{-i}}(Y_s^{\Gamma,j} -g_{ij}(s,X_s^{t,x})), \hspace{0.1cm} \forall s \leq T;}\\\\
 \textstyle {\int_0^T (Y_s^{\Gamma,i} -\displaystyle \max_{j\in \mathcal{I}^{-i}}(Y_s^{\Gamma,j} -g_{ij}(s,X_s^{t,x}))) dK_s^{\Gamma,i} = 0.}
\end{cases}
\end{equation}
First note that by Theorem \ref{existmono} the solution of this system \eqref{eq12} exists and is unique since the generators $(\bar f_i:=f_i(r,X_r^{t,x},\vec{\Gamma}_r,z))_{i\in \mi}$, which do not depend on $\vec y$, and the functions $(h_i)_{\ii}$ and $(g_{ij})_{i,j\in \mi}$ satisfy the assumptions (H1)-(H3) and (H4)-(i) as well. It is connected with the optimal switching problem in the way which we will describe now.

Let $\delta:= (\theta_k, \alpha_k)_{k \geq 0}$ be an admissible strategy of switching, i.e., $(\theta_k)_{k \geq 0}$ is an increasing sequence of stopping times with values in $[0,T]$ such that $\mathbb{P}[\theta_k < T, \forall k \geq 0] = 0$ and $ \forall k \geq 0$, $\alpha_k$ is a random variable $\mathcal{F}_{\theta_{k}}$-measurable with values in $\mathcal{I}.$

 Next with the admissible strategy $\delta:= (\theta_k, \alpha_k)_{k \geq 0}$ is associated a switching cost process $(A_s^{\delta})_{s \leq T}$ defined by:
 \begin{equation}\label{defcoutswitch}
 A_s^{\delta} := \displaystyle \sum_{k \geq 1} g_{\alpha_{k-1} \alpha_{k}}(\theta_k, X_{\theta_{k}}^{t,x})\mathbf{1}_{\lbrace \theta_{k} \leq s \rbrace} \mbox{ for }s < T,\mbox{ and } A_T^{\delta}= \lim_{s\rightarrow T}A_s^{\delta}.
 \end{equation}
Note that $(A_s^{\delta})_{s \leq T}$ is an RCLL process. Now for any fixed $s\leq T$ and $i \in \mathcal{I}$, let us denote by $\mathcal{A}_s^i$ the following set of admissible strategies :
\begin{equation*}
\mathcal{A}_s^i:= \lbrace\delta:= (\theta_k, \alpha_k)_{k \geq 0} \ \mbox{admissible strategy such that}\ \theta_0=s, \alpha_0=i \ \mbox{and}\ \mathbb{E}[(A_T^{\delta})^2] < \infty \rbrace.
\end{equation*}
Next let $\delta:= (\theta_k, \alpha_k)_{k \geq 0} \in \mathcal{A}_s^i$ and let us define the pair of adapted processes $(P^{\delta}, N^{\delta}):= (P_s^{\delta}, N_s^{\delta})_{s \leq T}$ as follows:
 \begin{equation}\label{eq15}\left\{\begin{array}{l}
 P^\delta \mbox{ is RCLL and } \E[\sup_{s\le T}|P^\d_s|^2]<\infty \,;
 N^{\delta}\in \hdd ;\\\\
P_s^{\delta} = h^{\delta}(X_T^{t,x}) + \int_s^T f^{\delta}(r,X_r^{t,x},\vec{\Gamma}_r                                                                                                                                                                                                 ,N_r^{\delta})dr - \int_s^T N_r^{\delta}dB_r - A_T^{\delta} + A_s^{\delta},\,\forall \,s\leq T, \end{array}\right.
 \end{equation}
 where
 \begin{align}\label{eqfd}
 h^{\delta}(x):= \displaystyle \sum_{k \geq 0}h_{\alpha_{k}}(x) \mathbf{1}_{[\theta_{k} \leq T < \theta_{k+1})}\mbox{ and }
f^{\delta}(s,x,\vec y,z) :=  \displaystyle \sum_{k \geq 0}f_{\alpha_{k}}(s,x,\vec y,z)\mathbf{1}_{[\theta_{k} \leq s < \theta_{k+1})}.
 \end{align}
Those series contain only a finite many terms since $\d$ is admissible and then $\mathbb{P}[\theta_n<T,\forall n\geq 0]=0$.

Next by a change of variable, the existence of $(P^{\delta} - A^{\delta},N^\d)$ stems from the standard existence result of solutions of BSDEs by Pardoux-Peng \cite{pardouxpeng} since its generator $ f^{\delta}(s,X_s^{t,x},\vec{\Gamma}_s                                                                                                                                                     ,z)$ is Lipschitz $\wt$ $z$ and $A^\d_T$ is square integrable. Then the solution of \eqref{eq15} follows.

We then have the following link between the solution of \eqref{eq12} and the value function of the optimal switching problem (see e.g. \cite{dhmz, hamadene2013viscosity, hu2010multi} for more details on this representation) :
 \begin{equation}\label{eq13}
 Y_s^{\Gamma,i} = \displaystyle \mbox{esssup}_{\delta \in \mathcal{A}_s^i}(P_s^{\delta} - A_s^{\delta}) = P_s^{\delta^{*}} - A_s^{\delta^{*}},
 \end{equation}
for some $\delta^{*} \in \mathcal{A}_s^i$, which means that $\d^*$ is an optimal strategy of the switching control problem.

Let us now show that $\Theta$ is a contraction in $\mathcal{H}^{2,m}$ with an appropriate norm. First let us introduce the following equivalent norm on this latter space: for any $\vec y\in \mathcal{H}^{2,m}$,
\begin{equation*}
\|\vec y\|_{\alpha}^2 := \mathbb{E}\Big[\int_0^T e^{\alpha s} |\vec y_s|^2ds\Big].
\end{equation*}
Next let ${}^1\vec \Gamma $, ${}^2\vec \Gamma$ be two processes of $\mathcal{H}^{2,m}$.
For $i\in \mi$, let us set: $$F_i(s,\omega, X_s^{t,x}(\omega),z):= f_i(s,X_s^{t,x}(\omega), {}^1 \vec \Gamma_s(\omega),z) \vee f_i(s,X_s^{t,x}(\omega), {}^2\vec \Gamma_s(\omega),z).$$
Since  $F_i$ satisfies (H1)-(H3) and (H4)-(i), we denote by $(\tilde{Y}^i, \tilde{Z}^i, \tilde{K}^i)_{i \in \mathcal{I}}$  the solution of the obliquely reflected BSDEs associated with $((F_i)_{i\in \mathcal{I}}, (g_{ij})_{i,j \in \mathcal{I}},(h_i)_{i \in \mathcal{I}})$. Moreover, once more, we have the following representation : $\forall s \leq T$,
\begin{equation}\label{eq14}
 \tilde{Y}_s^{i} = \displaystyle \ssp_{\delta \in \mathcal{A}_s^i}(\tilde{P}_s^{\delta} - A_s^{\delta})=(\tilde{P}_s^{{\tilde \delta}^*} - A_s^{{\tilde \delta}^*}),
\end{equation}
where for any strategy $\d$ of $\mathcal{A}_s^i$, the pair of processes $(\tilde{P}^{\delta}, \tilde N ^\d)$ is adapted and verifies:
\begin{equation}\label{eq15x}\left\{\begin{array}{l}
\tilde P^\delta \mbox{ is RCLL and } \E[\sup_{s\le T}|\tilde P^\d_s|^2]<\infty \,;
\tilde  N^{\delta}\in \hdd ;\\\\
\tilde P_s^{\delta} = h^{\delta}(X_T^{t,x}) + \int_s^T F^{\delta}(r,X_r^{t,x},\tilde N_r^{\delta})dr - \int_s^T \tilde N_r^{\delta}dB_r - A_T^{\delta} + A_s^{\delta},\,\forall s\leq T.
 \end{array}\right.
 \end{equation}
Here the generator $F^\d(...)$, associated with $\d$ is defined in the same way as in \eqref{eqfd} where we defined $f^\d(...)$ from $(f_i)_{\ii}$. Then by comparison (see Theorem\ref{existmono}), we have  $\forall s\leq T$, $Y_s^{{}^1 \Gamma,i} \leq \tilde{Y}_s^{i} $ and $Y_s^{{}^2 \Gamma,i} \leq \tilde{Y}_s^{i} $. This combined with \eqref{eq13} and \eqref{eq14}, leads to
$$
{}^1 \!P_s^{ {\tilde \delta}^*} - A_s^{{\tilde \delta}^*}\leq Y_s^{{}^1  \Gamma,i} \leq \tilde{P}_s^{{\tilde \delta}^*}- A_s^{{\tilde \delta}^*} \hspace{0.1cm} \mbox{ and }{}^2 \!P_s^{ {\tilde \delta}^*} - A_s^{{\tilde \delta}^*} \leq Y_s^{{}^2  \Gamma,i} \leq \tilde{P}_s^{{\tilde \delta}^*}- A_s^{{\tilde \delta}^*},$$
where ${}^1 \!P^{ {\tilde \delta}^*}$ (resp. ${}^2 \!P^{ {\tilde \delta}^*}$) is the first component the solution of the BSDE \eqref{eq15} with generator \\ $f^{ {\tilde \delta}^*}(s,X_s^{t,x},{}^1 \Gamma_s,z)$ (resp. $f^{ {\tilde \delta}^*}(s,X_s^{t,x},{}^2 \Gamma_s,z)$). We then have $$
|Y_s^{{}^1  \Gamma,i} - Y_s^{{}^2  \Gamma,i}| \leq |\tilde{P}_s^{{\tilde \delta}^*} - {}^1 \!P_s^{{\tilde \delta}^*}| + |\tilde{P}_s^{{\tilde \delta}^*} - {}^2 \!P_s^{\delta^{*}}|$$ and then
\begin{equation}\label{eq42}|Y_s^{{}^1  \Gamma,i} - Y_s^{{}^2  \Gamma,i}|^2 \leq 2\{|\tilde{P}_s^{{\tilde \delta}^*} - {}^1 \!P_s^{{\tilde \delta}^*}|^2 + |\tilde{P}_s^{{\tilde \delta}^*} - {}^2 \!P_s^{\delta^{*}}|^2\}.
\end{equation}
Next multiplying both members of the last inequality by $e^{\alpha s}$, using It\^o's formula (with the right hand-side) and the inequality $|x\vee y-x|\le |x-y|$, $\forall x,y\in \R$, to deduce that:
\begin{equation}\label{ineq:important}\begin{array}{l}\forall s\leq T,\,\,
\E[e^{\alpha s}|Y_s^{{}^1  \Gamma,i} - Y_s^{{}^2  \Gamma,i}|^2]\leq \frac{2C}{\alpha}\E[\int_s^Te^{\alpha r}|
{}^1  \Gamma_r-{}^2  \Gamma_r|^2dr]\end{array}
\end{equation}
where $C$ is a common Lipschitz constant of the $f_i's$ w.r.t $(\vec y,z)$ and $\alpha \geq C$. Next since \eqref{ineq:important} is valid for any $(s,i)\in [0,T]\times \mi$, we obtain by integration:
\begin{equation}\label{eq16}
\|\Theta({}^1  \Gamma ) - \Theta({}^2  \Gamma )\|^2_{\alpha} \leq \frac{2CTm}{\alpha} \|{}^1
\Gamma - {}^2  \Gamma\|^2_{\alpha}.
\end{equation}
Henceforth there exists some appropriate constant $\alpha_0 >0$ (it is enough to take
$\alpha_0=4CTm$) such that $\Theta$ is contraction on the Banach space $(\mathcal{H}^{2,m},\|.\|_{\alpha_0})$. Thus it has a fixed point which provides the unique solution of system \eqref{eqq2}.\qed

We next provide some properties of the solution of system (\ref{eqq2}) which will be useful later.

\begin{proposition} \label {defui} Assume (H1)-(H3). Then:
\medskip

\noindent i) There exist deterministic functions $(u^i)_{\ii}$ of polynomial growth, defined on $\spa$, such that:
$$
\forall \,\ii, Y^{i,t,x}_s=u^i(s,\xtx_s),\,\, ds\times d\mathbb{P}\mbox{ on }[t,T]\times \R^k.$$
ii) Assume moreover that $f_i(t,x,0,0)$ and $h_i(x)$ are bounded. Then the processes $Y^{i,t,x}$ and functions $u^i$, $\ii$, are also bounded.
\end{proposition}
{\bf \sol{Proof}}: First let us focus on the first point. Let $(\bar Y,\bar Z)$ be the solution of the following standard BSDE:
$$\left\{
\begin{array}{l}
\bar Y\in \ss, \bar Z\in \hdd;\\
\bar Y_s=\Phi (X^{t,x}_T)+\int_s^T
\Psi(r,X^{t,x}_r,\bar Y_r,\bar Z_r)dr-\int_s^T\bar Z_rdB_r, \,\forall \,s\leq T
\end{array}
\right.
$$where for any $(s,x,y,z)\in [0,T]\times \R^{k+1+d}$, $\Psi(s,x,y,z):=\bar Cm|y|+\bar C|z|+\sum_{i=1,m}|f_i(s,x,0,\dots,0,)|$ and \\$\Phi(x):=\sum_{i=1,m}|h_i(x)|$. The constant $\bar C:=C_1+...+C_m$ with, for any $\ii$, $C_i$ is the Lipschitz constant of $f_i$ $\wt$ $(\vec y,z)$.

First note that since $\Psi \ge 0$ and $\Phi\ge 0$ then $\bar Y\ge 0$. Next as we are in the Markovian framework of randomness and since $\Phi$ and $\Psi(t,x,0,0)$ are of polynomial growth, then there exists a deterministic function $v(t,x)$ of polynomial growth (see e.g. \cite{kpq}) such that:
$$
\forall \, s\in [t,T], \bar Y^{t,x}_s=v(s,\xtx_s).
$$
Next let us set, for $i\in \mi$,
$$
\underbar Y_i=\bar Y, \underbar Z_i=\bar Z \mbox{ and } \underbar K_i=0.
$$Therefore, since $g_{ij}\ge 0$ for any $i,j\in \mi$, $(\underbar Y_i,\underbar Z_i,\underbar K_i)_{i\in \mi}$ is a solution of the following system: for any $\ii$ and $s\le T$,
\begin{equation}\label{eqbaryi}\left\{
\begin{array}{l}
\underbar Y_i(s)=\Phi(X^{t,x}_T)+\int_s^T
\Psi(r,X^{t,x}_r,\underbar Y_i(r),\underbar Z_i(r))dr+\underbar K_i(T)-\underbar K_i(s)-\int_s^T\underbar Z_i(r)B_r\,\,;\\\\
\underbar Y_i(s)\geq \max_{j\neq i}\{\underbar Y_j(s)-g_{ij}(s)\};\\\\
\int_0^T(\underbar Y_i(s)-\max_{j\neq i}\{\underbar Y_j(s)-g_{ij}(s)\})d\underbar K_i(s)=0.
\end{array}
\right.
\end{equation}In the same way let us set for any $i\in \mi$,
$$
\hat Y_i=-\bar Y, \hat Z_i=-\bar Z \mbox{ and } \hat K_i=0,$$ then $(\hat Y_i,\hat Z_i,\hat K_i)_{i\in \mi}$ is a solution of the following system: for any $\ii$ and $s\le T$,
$$\left\{
\begin{array}{l}
\hat Y_i(s)=-\Phi(X^{t,x}_T)-\int_s^T\Psi(r,X^{t,x}_r,-\bar Y_r,-\hat Z_r)dr+\hat K_i(T)-\hat K_i(s)-\int_s^T\hat Z_i(r)B_r\,\,;\\\\
\hat Y_i(s)\geq \max_{j\neq i}\{\hat Y_j(s)-g_{ij}(s)\};\\\\
\int_0^T(\hat Y_i(s)-\max_{j\neq i}\{\hat Y_j(s)-g_{ij}(s)\})d\hat K_i(s)=0.
\end{array}
\right.$$
Next let us consider the following sequence of processes $((\tilde Y_k^i,\tilde Z_k^i,\tilde K_k^i)_{\ii})_{k\ge 0}$:
$$
\tilde Y_0^i=0 \mbox{ for all }\ii \mbox{ and for }k\geq 1,\,\,(\tilde Y_k^i)_{\ii}=\Theta((\tilde Y_{k-1}^i)_{\ii})
$$
where $\Theta$ is the mapping defined in \eqref{defteta}. Therefore, as pointed out in the proof of Theorem \ref{thmexistence}, the sequence $((\tilde Y_k^i,\tilde Z_k^i,\tilde K_k^i)_{\ii})_{k\ge 0}$ converges to
$(Y^i,Z^i,K^i)_{\ii}$ in $(\mathcal{H}^{2,m},\|.\|_{\alpha_0})$ since $\Theta$ is a contraction in this latter complete normed space. On the other hand by an induction argument on $k$ and by using the comparison result of Theorem \ref{existmono}-ii), we have that:
\begin{equation}\lb {recurrence}\forall k\geq 0, \,\forall \ii,\,\,-\bar Y^i=\hat Y\leq \tilde Y_k^i\leq \underbar Y_i=\bar Y.
\end{equation}
Indeed for $k=0$, this obviously holds since $\bar Y\geq 0$. Next suppose that \eqref{recurrence} holds for some $k-1$ with $k\geq 1$. Then by a linearization procedure of $f_i$,  which is possible since it is Lipschitz $\wt$ $(\vec y,z)$, we have: for any $\ii$,
$$f_i(s,\xtx_s,(\tilde Y_{k-1}^i(s))_{\ii},z)=
f_i(s,\xtx_s,0,0)+\sum_{l=1,m}a_s^{k,i,l}\tilde Y_{k-1}^l (s)+ b_s^{k,i,l}z$$
where $a^{k,i,l}\in \R$ and $b^{k,i,l}\in \R^d$ are $\cal P$-measurable processes, bounded by the Lipschitz constant of $f_i$. Therefore, using the induction hypothesis, we obtain: $$|f_i(s,\xtx_s,(\tilde Y_{k-1}^i(s))_{\ii},z)|\le \Psi(s,\xtx_s,\bar Y_s,z).$$
Finally by the comparison argument of Theorem \ref{existmono}-ii) (see also \cite{hamadene2010switching}, Cor. 3.4, pp.411), we get: \\$\forall \ii,\, \tilde Y_k^i\leq \underbar Y'_i$ where $(\underbar Y'_i,\underbar Z'_i,\underbar K'_i)_{i\in \mi}$ is the unique solution of the system of type (\ref{eqq2}) associated with $((\underbar f_i=\Psi(s,\xtx_s,\underbar Y_s,z))_{\ii},(\underbar h_i=\Phi(x))_{\ii}, (g_{ij}(s,\xtx_s))_{i,j\in \mi})$. But the solution of this latter system is unique (Theorem \ref{thmexistence}) and by \eqref{eqbaryi}, $(\underbar Y_i,\underbar Z_i,\underbar K_i)_{i\in \mi}$ is also a solution. Therefore for any $\ii$, $\underbar Y'_i=\underbar Y_i$ and then
 $\forall \ii,\, \tilde Y_k^i\leq \underbar Y_i=\bar Y$. In the same way one can show that $\forall \ii,\, \tilde Y_k^i\geq \hat Y_i=-\bar Y$. Therefore \eqref{recurrence} holds true for any $k\ge 0$.

Next, once more, since we are in the Markovian framework of randomness, and using an induction argument on $k$ we deduce the existence of deterministic continuous functions of polynomial growth $u^{i,k}(t,x)$ (see e.g. \cite{hamadene2013viscosity}, Cor.2, pp.182), $\ii$, such that for any $\ii$, $\tx$,
\begin{equation}\lb{exp:uik}
\tilde Y^{i}_k(s)=u^{i,k}(s,\xtx_s), \,\,\forall s\in [t,T].
\end{equation}
By \eqref{recurrence}, in taking $s=t$, we obtain: for any $k\ge 0$, $\ii$ and $\tx$,
\begin{equation}\lb{estim:uik}
|u^{i,k}(t,x)|\leq v(t,x).
\end{equation}
Next by using the inequality \eqref{ineq:important} at $s=t$ we deduce that for any $\ii$, $k,p\geq 1$
 \begin{equation}\label{ineq:important2}\begin{array}{ll}|u^{i,k}(t,x) - u^{i,p}(t,x)|^2] &\leq\frac{2C}{\alpha_0}\E[\int_t^Te^{\alpha_0 (r-t)}\sum_{j=1,m}|\tilde Y^{j}_{k-1}(r) - \tilde Y^{j}_{p-1}(r)|^2dr].\\\\
 &\le  \frac{2C}{\alpha_0}\E[\int_t^Te^{\alpha_0 (r-t)}\sum_{j=1,m}|u^{j,k-1}(r,\xtx_r) - u^{j,p-1}(r,\xtx_r)|^2dr].\end{array}
\end{equation}
As $((\tilde Y_k^i)_{\ii})_k$ is a Cauchy sequence in $(\mathcal{H}^{2,m},\|.\|_{\alpha_0})$, then $((u^{i,k})_{\ii})_k$ is a Cauchy sequence pointwisely. This implies the existence of deterministic functions $(u^i)_{\ii}$ such that for any $\ii$ and $\tx$, $u^{i,k}(t,x)$ converges in $k$ to $u^{i}(t,x)$. Moreover by
\eqref{estim:uik}, $u^i$ is of polynomial growth since $v$ is so and finally by \eqref{exp:uik}, $Y^{i,t,x}_s=u^i(s,\xtx_s)$, $ds\times d\mathbb{P}$ on $\spat$.
\medskip

We now deal with the second point. Assume that $f_i(t,x,0,0)$ and $h_i(x)$ are bounded. Then the solution $\bar Y$ is bounded. This is obtained by a change of probability and by multiplying both hand-sides of the equation by $e^{-m\bar Cs}$, conditionning and taking into account of $\bar Y\geq 0$. Therefore the deterministic function $v$ is a also bounded. Consequently, $u^{i,k}$ are uniformly bounded and so are $u^i$, $\ii$.\qed
\begin{remarque}At this point we do not know whether the functions $u^i$, $\ii$, are continuous or not. However we will show later that they can be chosen continuous.\qed
\end{remarque}

\section{The main result : Existence and uniqueness of the viscosity solution for system of PDEs with interconnected obstacles}

In this section, we study the existence and uniqueness in viscosity sense of the solution of  the system of $m$ partial differential equations with interconnected obstacles \eqref{eq12}. The candidate to be the solution are the functions $(u^1,\dots,u^m)$ defined in Proposition \ref{defui} by which we represent $(Y^i)_{\ii}$.  So, firstly we are going to show that those functions $u^i$, $\ii$, can be chosen continuous.
\begin{proposition}Assume that (H1)-(H3) hold. Then we can choose the functions $u^i$, $\ii$, defined in Proposition \ref{defui}, continuous in $(t,x)$ and of polynomial growth.
\end{proposition}
\textbf{Proof}: It will be given in two steps. In the first one we are going to suppose moreover that $h_i$ and $f_i(t,x,0,0)$, $\ii$, are bounded. Later on we deal with the general case, i.e., without assuming the boundedness of those latter functions.

\noindent \tbf{\sol{Step 1}}: Suppose that for any $\ii$, $h_i$ and $f_i(t,x,0,0)$ are bounded.

\noindent Recall the continuous functions $u^{i,k}$, $\ii$ and $k\ge 0$, defined in \eqref{exp:uik}. By \eqref{ineq:important2} they verify: $\forall k\geq 1$, $\ii$ and $\tx$,
\begin{equation}\label{ineq:important3}\begin{array}{l}|u^{i,k}(t,x) - u^{i,p}(t,x)|^2 \leq\frac{2C}{\alpha}\E[\int_t^Te^{\alpha (r-t)}\sum_{j=1,m}|u^{j,k-1}(r,\xtx_r) - u^{j,p-1}(r,\xtx_r)|^2dr]\end{array}
\end{equation}
where, as pointed out in the proof of Theorem \ref{thmexistence}, $C$ is the uniform Lipschitz constant of $f_i's$ $\wt$ $(\vec y,z)$ and $\alpha \geq C$.

On the other hand we know, by Proposition \ref{defui}-ii), that $u^{i,k}$ are uniformly
bounded for any $\ii$ and $k\geq 0$. Now let us take $\alpha=C$ and let $\eta$ be a constant such that $2C m (e^{C \eta }-1)= \frac{3}{4}$ and finally let us set
$$\|u^{i,k}- u^{i,p}\|_{\infty,\eta}:=\sup_{(t,x)\in [T-\eta, T]\times \R^k}|u^{i,k}(t,x) - u^{i,p}(t,x)|.
$$
Therefore, we deduce from \eqref{ineq:important2}, that for any $k,p\ge 1$,
$$
\sum_{i=1,m}\|u^{i,k}- u^{i,p}\|^2_{\infty,\eta}\leq  \frac{3}{4}\sum_{i=1,m}\|u^{i,k-1}- u^{i,p-1}\|^2_{\infty,\eta}
$$which means that the sequence $((u^{i,k})_{\ii})_{k\geq 0}$ is uniformly convergent in $[T-\d,T]\times \R^k$. Thus, their limits, i.e., the functions $(u^{i})_{\ii}$ are also continuous on the set $[T-\d,T]\times \R^k$.

Next let $s\in [T-2\eta, T-\eta]$, then once more by \eqref{ineq:important3} we have:
\begin{align}\label{ineq:important4}&|u^{i,k}(t,x) - u^{i,p}(t,x)|^2 \\ &\qquad \leq\frac{2C}{\alpha}\E[\int_t^{T-\eta}e^{\alpha (r-t)}\sum_{j=1,m}|u^{j,k-1}(s,\xtx_s) - u^{j,p-1}(s,\xtx_s)|^2dr]+\frac{3}{4}\sum_{i=1,m}\|u^{i,k-1}- u^{i,p-1}\|^2_{\infty,\eta}.\nonumber
\end{align}
And then if we set
$$\|u^{i,k}- u^{i,p}\|_{\infty,2\eta}:=\sup_{(t,x)\in [T-2\eta, T-\eta]\times \R^k}|u^{i,k}(t,x) - u^{i,p}(t,x)|,
$$we obtain:
$$
\sum_{i=1,m}\|u^{i,k}- u^{i,p}\|^2_{\infty,2\eta}\leq  \frac{3}{4}\sum_{i=1,m}\|u^{i,k-1}- u^{i,p-1}\|^2_{\infty,2\eta}+\frac{3}{4}\sum_{i=1,m}\|u^{i,k-1}- u^{i,p-1}\|^2_{\infty,\eta}.
$$
It implies that
$$
\limsup_{k,p\rightarrow \infty}\sum_{i=1,m}\|u^{i,k}- u^{i,p}\|^2_{\infty,2\eta}\leq  \frac{3}{4}\limsup_{k,p\rightarrow \infty}\sum_{i=1,m}\|u^{i,k-1}- u^{i,p-1}\|^2_{\infty,2\eta}$$
since $\limsup_{k,p\rightarrow \infty}\sum_{i=1,m}\|u^{i,k-1}- u^{i,p-1}\|^2_{\infty,\eta}=0.$ Therefore $$
\limsup_{k,p\rightarrow \infty}\sum_{i=1,m}\|u^{i,k}- u^{i,p}\|^2_{\infty,2\eta}=0.
$$Consequently the sequence $((u^{i,k})_{\ii})_{k\geq 0}$ is uniformly convergent in $[T-2\eta,T-\eta]\times \R^k$. Thus, their limits, the functions $(u^{i})_{\ii}$ are also continuous in $[T-2\eta ,T-\eta]\times \R^k$, which implies that $(u^{i})_{\ii}$ are continuous in $[T-2\eta ,T]\times \R^k$. Continuing now this reasoning as many times as necessary on
$[T-3\eta,T-2\eta]$, $[T-4\eta ,T-3\eta]$ etc. we obtain the continuity of $(u^{i})_{\ii}$ in $[0, T]\times \R^k$.

\nd  \textbf{Step 2 :} We now deal with the general case. Firstly by (H1)-iii), (H2) and (H3), there exist two constants $C$ and $p\in \N$ such
$f_i(t,x,0,...,0)$, $h_i(x)$ and $g_{ij}(t,x)$ are of polynomial growth, i.e., for any $(t,x) \in [0,T]\times \R^k$,
\begin{equation}\lb{cdtcroispoly}
|f_i(t,x,0,...,0)| +  |h_i(x)| + |g_{ij}(t,x)| \leq C(1+|x|^p).
\end{equation}
To proceed for $s\in [t,T]$ let us define,
\begin{equation*}
\overline{Y}_s^{i} := Y_s^{i} \varphi(X_s^{t,x}),
\end{equation*}
where for $ x\in \R$, $\varphi(x):= \frac{1}{(1+|x|^2)^p}$ ($p$ is the same constant as in \eqref{cdtcroispoly}). Then by the integration-by-parts formula we have:
\begin{align*}
d\overline{Y}_s^{i}  &= \varphi(X_s^{t,x})dY_s^{i}  + Y_s^{i} d\varphi(X_s^{t,x}) + d\langle Y^{i},\varphi(X^{t,x})\rangle_s \\\\
& = \varphi(X_s^{t,x}) \{ -f_i(s,X_s^{t,x},(Y_s^{k})_{k=1,...,m},Z_s^{i})ds - dK_s^{i} +  Z_s^{i}dB_s \}\\
  & \qquad \qquad + Y_s^{i} \{\mathcal{L}\varphi(X_s^{t,x})ds + D_x\varphi(X_s^{t,x})\sigma(s,X_s^{t,x})dB_s \} + Z_s^{i} D_x\varphi(X_s^{t,x}) \sigma(s,X_s^{t,x}) ds\\\\
& = \{ - \varphi(X_s^{t,x}) f_i(s,X_s^{t,x},(Y_s^{k})_{k=1,...,m},Z_s^{i,n})+  \mathcal{L}\varphi(X_s^{t,x}) Y_s^{i}+D_x\varphi(X_s^{t,x}) \sigma(s,X_s^{t,x})Z_s^{i}\}ds\\
&  \qquad \qquad -\varphi(X_s^{t,x}) dK_s^{i}  + \{ Z_s^{i}\varphi(X_s^{t,x}) + Y_s^{i} D_x\varphi(X_s^{t,x}) \sigma(s,X_s^{t,x}) \}dB_s,
\end{align*}
where $\mathcal{L}\varphi$ is given in \eqref{generinf}. Next let us set, for  $s\in [t, T]$,
\begin{align*}
d\overline{K}_s^{i} :=\varphi(X_s^{t,x})dK_s^{i} \mbox{ and }
\overline{Z}_s^{i} := Z_s^{i}\varphi(X_s^{t,x}) + Y_s^{i} D_x\varphi(X_s^{t,x}) \sigma(s,X_s^{t,x}).
\end{align*}
Then  $((\overline{Y}^{i} ,\overline{Z}^{i} ,\overline{K}^{i} ))_{\ii}$ satisfies: $\forall s\in [t,T]$,
\begin{equation}
\begin{cases}
\overline{Y}_s^{i} = \overline{h}_i(X_T^{t,x})+\int_s^T\overline{f}_i(r,X_r^{t,x},(\overline{Y}_r^{k})_{k=1,...,m},\overline{Z}_r^{i})dr +\overline{K}_T^{i}- \overline{K}_s^{i} -\int_s^T \overline{Z}_r^{i}dB_r, \\\\
\overline{Y}_s^{i} \geqslant  \displaystyle \max_{j \in \mathcal{I}^{-i}}(\overline{Y}_s^{j} -\overline{g}_{ij}(s,X_s^{t,x})),\\\\
\int_t^T (\overline{Y}_s^{i} -\displaystyle \max_{j\in \mathcal{I}^{-i}}(\overline{Y}_s^{j} -\overline{ g}_{ij}(s,X_s^{t,x}))) d\overline{K}_s^{i} = 0,
\end{cases}
\end{equation}
where for any $i,j\in \mi$,
$$
\overline{h}_i(X_T^{t,x}):= h_i(X_T^{t,x})\varphi(X_T^{t,x}),\,\, \overline{g}_{ij}(s,X_s^{t,x}) := g_{ij}(s,X_s^{t,x})\varphi(X_s^{t,x}), $$
and
\begin{align*}
&\overline{f}_i (s,x,\vec{y}, z):= \varphi(x)f_i (s,x,\varphi^{-1}(x)\vec y, \varphi^{-1}(x)z-D_x\varphi(x) \sigma(s,x)\varphi^{-1}(x)y^i) )\\
& \qquad \qquad \qquad \qquad+ \mathcal{L}\varphi(x)\varphi^{-1}(x)y^i+  D_x\varphi(x) \sigma(s,x)\varphi^{-1}(x)\{z- D_x\varphi(x) \sigma(s,x)\varphi^{-1}(x)y^i\}.
\end{align*}
Here let us notice that the functions $\overline{f}_i(t,x,0,0)$, $\overline{g}_{ij}$ and $\overline{h}_i$ are bounded. Then by the result of the first step, there exists bounded continuous functions $(\bar u^i)_{\ii}$ such that for any $\tx$, and $s\in [t,T]$, $\bar Y^i_s=\bar u^i(s,\xtx_s)$, $\forall i\in \mi$.  Thus
for any $\tx$, and $s\in [t,T]$, $Y^i_s=\varphi^{-1}(\xtx_s)\bar u^i(s,\xtx_s)$, $\forall i\in \mi$. Then it is enough to take $u^i(t,x):=\varphi^{-1}(x)\bar u^i(t,x)$, $\tx$ and $\ii$, which are continuous functions and of polynomial growth.
\qed

We are now ready to give the main result of this paper. Let $(Y^i,Z^i,K^i)_{\ii}$ be the unique solution of \eqref{eqq2} and let $(u^i)_{\ii}$ be the continuous functions with polynomial growth such that for any $\tx$, $\ii$ and $s\in [t,T]$, $$Y^{i,t,x}_s=u^i(s,\xtx_s).$$
We then have:
\begin{theoreme} The functions $(u^i)_{\ii}$ is a solution in viscosity sense of system \eqref{eqq12}. Moreover it is unique in the class of continuous functions of polynomial growth.
\end{theoreme}
\sol{\tbf{Proof}}: First let us show that $(u^i)_{\ii}$ is a viscosity solution of system
\eqref{eqq12}.

Recall that  $(Y^i,Z^i,K^i)_{\ii}$ is a solution of the system of reflected BSDEs with interconnected obstacles \eqref{eqq2} and for any $\tx$, $\ii$ and $s\in [t,T]$, $Y^{i,t,x}_s=u^i(s,\xtx_s)$. Then $(Y^i,Z^i,K^i)_{\ii}$ verify: for any $s\in [t,T]$ and $\ii$,
\begin{equation}\label{eqq2u}
\begin{cases}
\textstyle {Y_s^{i,t,x} = h_i(X_T^{t,x}) + \int_s^T f_i(r,X_r^{t,x},(u^k(r,\xtx_r))_{k\in \mi},Z_r^{i,t,x})dr + K_T^{i,t,x} - K_s^{i,t,x}- \int_s^T Z_r^{i,t,x}dB_r,}\\\\
 \textstyle {Y_s^{i,t,x} \geqslant  \displaystyle \max_{j \in \mathcal{I}^{-i}}(u^j(s,\xtx_s) -g_{ij}(s,X_s^{t,x}))}\mbox{ and }\textstyle {\int_t^T (Y_r^{i,t,x} -\displaystyle \max_{j\in \mathcal{I}^{-i}}(u^j(r,\xtx_r)  -g_{ij}(r,X_r^{t,x}))) dK_r^{i,t,x} = 0.}
\end{cases}
\end{equation}
But system \eqref{eqq2u} is decoupled and using a result by El-Karoui et al. (Theorem 8.5 in \cite{kkppq}) one obtains that, for any $i_0$, $u^{i_0}$ is a solution in viscosity sense of the following PDE with obstacle:
\begin{equation}\label{eq7}
\begin{cases}
\min \lbrace u^{i_0}(t,x) - \displaystyle \max_{j \in {\mathcal{I}^{-i}}}(u^j(t,x)-g_{ij}(t,x)) ;\\
 \quad \quad -\partial_tu^{i_0}(t,x) - \mathcal{L}u^{i_0}(t,x) - f_{i_0}(t,x,({u}^k(t,x))_{k=1,,...,m},(\sigma^\top D_xu^{i_0})(t,x))\rbrace = 0  ;\\\\
u^{i_0}(T,x) = h_i(x).
\end{cases}
\end{equation}
As $i_0$ is arbitrary in $\mi$, then the functions $(u^i)_{\ii}$ is a solution in viscosity sense of \eqref{eqq12}.

Next let us show that $(u^i)_{\ii}$  is the unique solution in the class of continuous functions with polynomial growth. It is based on the uniqueness of the solution of the system of reflected BSDEs with interconnected obstacles \eqref{eqq2}.

So suppose that there exists another continuous with polynomial growth solution $(\tilde{u}^i)_{i=1,...,m}$ of \eqref{eqq12}, i.e.,
for any $i \in \mathcal{I}$,
\begin{equation}\label{eqq13}
\begin{cases}
\min \lbrace \tilde u^i(t,x) - \displaystyle \max_{j \in {\mathcal{I}^{-i}}}(\tilde u^j(t,x)-g_{ij}(t,x)) ;\\
 \quad \quad -\partial_t\tilde u^i(t,x) - \mathcal{L}\tilde u^i(t,x) - f_i(t,x,(\tilde u^k(t,x))_{k=1,...,m},(\sigma^\top D_x\tilde u^i)(t,x))\rbrace = 0 \,\, ;\\
\tilde u^i(T,x) = h_i(x).
\end{cases}
\end{equation}
Let $(\tilde{Y}^{i})_{i \in \mathcal{I}} \in \mathcal{H}^{2,m}$ be such that for any $i \in \mathcal{I}$ and $s\in [t,T]$,
\begin{equation*}
\tilde{Y}_s^{i,t,x} = \tilde{u}^i(s,X_s^{t,x}).
\end{equation*}
Next let us define $(\overline{Y}^{i,t,x})_{i \in \mathcal{I}}$ as follows:
\begin{equation}\lb{eqrep13}
(\overline{Y}^{i,t,x})_{i \in \mathcal{I}}= \Theta((\tilde{Y}_s^{i,t,x})_{i \in \mathcal{I}}),
\end{equation}
 that is to say, $(\overline{Y}^{i,t,x}, \overline{Z}^{i,t,x},\overline{K}^{i,t,x})_{i \in \mathcal{I}}$ is the solution of the following system of reflected BSDEs with oblique reflection: $\forall s\in [t,T]$,
 \begin{equation}\label{eqq22xx}
\begin{cases}
\textstyle {\overline{Y}_s^{i,t,x} = h_i(X_T^{t,x}) + \int_s^T f_i(r,X_r^{t,x},(\tilde{u}^k(s,X_s^{t,x})_{k=1,...m},\overline{Z}_r^{i,t,x})dr + \overline{K}_T^{i,t,x} - \overline{K}_s^{i,t,x}- \int_s^T \overline{Z}_r^{i,t,x}dB_r}\,\,;\\\\
 \textstyle {\overline{Y}_s^{i,t,x} \geqslant  \displaystyle \max_{j \in \mathcal{I}^{-i}}(\overline{Y}_s^{j,t,x} -g_{ij}(s,X_s^{t,x}))}\mbox{ and }\textstyle {\int_t^T (\overline{Y}_s^{i,t,x} -\displaystyle \max_{j\in \mathcal{I}^{-i}}(\overline{Y}_s^{j,t,x} -g_{ij}(s,X_s^{t,x}))) d\overline{K}_s^{i,t,x} = 0.}
\end{cases}
\end{equation}
As the deterministic functions $(\tilde{u}^i)_{i=1,...,m}$ are continuous and of polynomial growth, then by using a result by Hamad\`ene-Morlais (\cite{hamadene2013viscosity}, Theorem 1), one can infer the existence of deterministic continuous functions with polynomial growth $(v^i)_{i=1,...,m}$ such that: $\forall i \in \mathcal{I}$ and $ s\in [t,T]$,
\begin{equation*}
\overline{Y}_s^{i,t,x} = v^i(s,X_s^{t,x}).
\end{equation*}
Moreover, $(v^i)_{i=1,...,m}$ is the unique viscosity solution (in the class of functions with polynomial growth) of the following system of PDEs with interconnected obstacles : $\forall i = {1,...,m}$
\begin{equation}\label{eq7}
\begin{cases}
\min \lbrace v^i(t,x) - \displaystyle \max_{j \in {\mathcal{I}^{-i}}}(v^j(t,x)-g_{ij}(t,x)) ;\\
 \quad \quad -\partial_tv^i(t,x) - \mathcal{L}v^i(t,x) - f_i(t,x,(\tilde{u}^k(t,x))_{k=1,,...,m},(\sigma^\top D_xv^i)(t,x))\rbrace = 0  \,\,;\\\\
v^i(T,x) = h_i(x).
\end{cases}
\end{equation}
Let us notice that, in system (\ref{eq7}), in the arguments of $f_i$ we have $\tilde{u}^k$ and not $v^k$. Now as the functions $(\tilde{u}^i)_{i=1,...,m}$ solve system (\ref{eq7}), hence by uniqueness of the solution of this system
\eqref{eq7} (see \cite{hamadene2013viscosity}, Thm. 1, pp.175), one deduces that
\begin{equation*}
\tilde{u}^i = v^i \,\, \mbox{ and then }\tilde{Y}^{i,t,x}=\overline{Y}^{i,t,x}, \,\,\forall \ii .
\end{equation*}
Therefore $(\tilde{Y}_s^{i,t,x})_{i \in \mathcal{I}}$ verify
\begin{equation*}
(\tilde{Y}^{i,t,x})_{i \in \mathcal{I}}= \Theta((\tilde{Y}_s^{i,t,x})_{i \in \mathcal{I}}).
\end{equation*}
But $(Y^i)_{i \in \mathcal{I}}$ is the unique fixed point of $\Theta$ in $(\mathcal{H}^{2,m}, \|.\|_{\alpha_0})$ then we have that for any $s \in [t,T]$ and $ i\in \mathcal{I}$, $
\tilde{Y}_s^{i,t,x} = Y_s^{i,t,x}.$ Henceforth, in taking $s=t$, we obtain that for any $i \in \mathcal{I}$ and
$\tx$, $\tilde{u}^i (t,x)= u^i(t,x)$. Thus $({u}^i)_{i=1,...,m}$ is the unique solution of system (\ref{eq7}) in the class of continuous functions with polynomial growth. \qed

\end{document}